\documentclass[12pt,oneside]{article}
\usepackage{amssymb,amsmath}
\usepackage{url}
\usepackage[english]{babel}
\renewcommand{\le}{\leqslant}
\renewcommand{\ge}{\geqslant}
\newcommand{\RR}{\mathbb{R}}
\newcommand{\ZZ}{\mathbb{Z}}

\newcommand{\NN}{\mathbb{N}}

\newcommand{\va}{\mathbf{a}}
\newcommand{\vb}{\mathbf{b}}
\newcommand{\vf}{\mathbf{f}}

\newcommand{\HHH}{\mathcal{H}}

\newcommand{\PPP}{\mathcal{P}}
\newcommand{\FFF}{\mathcal{F}}

\newcommand{\CCC}{\mathcal{C}}

\newtheorem{theorem}{Theorem}
\newtheorem{lemma}{Lemma}
\newtheorem{definition}{Definition}
\newtheorem{prop}{Proposition}
\newtheorem{corl}{Corollary}
\begin{document}

\title{Inhomogeneous Diophantine approximation on curves and Hausdorff dimension}

\author{Dzmitry Badziahin\footnote{Supported by EPSRC Grant EP/E061613/1} \\ {\small\sc York} }

\date{}

\maketitle

{\abstract{\scriptsize The goal of  this paper is to develop a
coherent theory for inhomogeneous Diophantine approximation on
curves in $\RR^n$  akin to the well established homogeneous theory.
More specifically, the measure theoretic results obtained generalize
the fundamental homogeneous theorems of R.C. Baker (1978), Dodson,
Dickinson (2000) and Beresnevich, Bernik, Kleinbock, Margulis
(2002). In the case of planar curves,  the complete Hausdorff
dimension theory is developed.}}

\section{Introduction}
Throughout $\psi:\RR^+\to \RR^+:= (0,+\infty) $ denotes a decreasing
function and will be referred to as an \emph{approximation
function}. Let $\vf=(f_1,\dots,f_n):I\to \RR^n$ be a $C^{(n)}$ map
defined on an interval $I\subset\RR$ and $\lambda:I\to\RR$ be a
function. For reasons that will soon be apparent, the function
$\lambda$ will be referred to as an \emph{inhomogeneous function}.
Let $A_n(\psi,\lambda)$ be the set of $x\in I$ such that the
inequality
\begin{equation}\label{e:001}
\|\va\cdot\vf(x)+\lambda(x)\| < \psi(|\va|)
\end{equation}
holds for infinitely many $\va\in\ZZ^n\setminus\{0\}$, where
$\|\cdot\|$ denotes the distance to the nearest integer,
$|\va|:=\max\{|a_1|,\ldots,|a_n|\}$ and $\va\cdot\vb$ stands for the
standard inner product of vectors $\va$ and $\vb$ in $\ZZ^n$. In the
special case when $\psi(h)=\psi_v(h):=h^{-v}$ for some fixed
positive $v$ we will denote $A_n(\psi_v,\lambda)$ by
$A_n(v,\lambda)$. Furthermore, in the case when the inhomogeneous
function  $\lambda$ is identically zero we write $A_n(\psi)$ for
$A_n(\psi,\lambda)$ and $A_n(v)$ for $A_n(v,\lambda)$.

By definition, $A_n(\psi,\lambda)$  is the set of $x\in I$ such that
the corresponding point $\vf(x)$ lying on the curve
\begin{equation}\label{e:002}
\CCC:=\big\{(f_1(x),f_2(x),\ldots,f_n(x)):x\in I\big\}\subset \RR^n
\end{equation}
satisfies the Diophantine condition arising from (\ref{e:001}).
Within the homogeneous setup ($\lambda\equiv0$), investigating the
measure theoretic properties of $A_n(\psi)$ dates back to 1932 and a
famous problem of Mahler \cite{Mahler}. The problem states that when
(\ref{e:002}) is the Veronese curve given by $(x,x^2,\dots,x^n)$
then $A_n(v)$ is of Lebesgue measure zero whenever $v>n$. Mahler's
problem was eventually settled by Sprindzuk \cite{sprindz2} in 1965
and subsequently Schmidt \cite{S} extended the result to the case of
arbitrary planar curves (\emph{i.e.} $n=2$) with non-vanishing
curvature. These two major results led to what is currently known as
the (homogeneous)  theory of Diophantine approximation on manifolds
\cite{BD}.

Diophantine approximation on manifolds has been an extremely active
research area over the past 10 years or so.   Rather than describe
the activity in detail, we refer the reader to research articles
\cite{beresnevich-alg, beresnevich, BBKM, BDV07, BKM, KM,
Vaughan-Velani-2007} and the surveys
\cite{Beresnevich-Bernik-Dodson-02:MR1975457, Kleinbock-survey,
Margulis}. Nevertheless, it is worth singling out the pioneering
work of Kleinbock and Margulis \cite{KM} in which the fundamental
Baker-Sprindzuk conjecture is established. This has undoubtedly
acted as the catalyst to the works cited above which together
constitute a  coherent homogeneous theory for Diophantine
approximation on manifolds. The situation for the inhomogeneous
theory is quite different.  Indeed,  the  inhomogeneous analogue  of
the Baker-Sprindzuk conjecture \cite{Beresnevich-Velani-Moscow,
Beresnevich-Velani-08-Inhom} has only just been established in 2008.

The aim of this paper is to develop a coherent theory for
inhomogeneous Diophantine approximation on curves akin to the well
established homogeneous theory. More precisely, Hausdorff measure
theoretic statements for the sets $A_n(\psi,\lambda)$ are obtained.
In particular, a complete metric theory is established in the case
of planar curves ($n=2$). In short, the results constitute the first
precise and general statements in the theory of inhomogeneous
Diophantine approximation on manifolds.

\subsection{Main results and corollaries}

Before we proceed with the statement of the results, we introduce
some useful notation and recall some standard definitions.

\noindent The curve $\mathcal{C}$ given by (\ref{e:002}) is called
\emph{non-degenerate} at $x\in I$ if the Wronskian
$$w(f'_1,\dots,f'_n)(x):=\det (f^{(j)}_i(x))_{1\le i,j\le n}$$  does
not vanish. We say that $\mathcal{C}$ is \emph{non-degenerate}\/ if
it is non-degenerate at almost every point $x\in I$. Given a set
$X\subset\RR^n$ and a real number  $s>0$, $\HHH^s(X)$ will denote
the $s$-dimensional Hausdorff measure of $X$ and $\dim X$ will
denote the Hausdorff dimension of $X$. The latter is defined to be
the infimum over $s$ such that $\HHH^s(X)=0$. For the formal
definitions and properties of Hausdorff measure and dimension see
\cite{Falconer}.

\bigskip

\noindent\textbf{Lower bounds.} Our first result enables us to
deduce lower bounds for $\dim A_n(v,\lambda)$ and represents an
inhomogeneous version of the homogeneous theorem established by
Dodson and Dickinson \cite{dod_dick}. Furthermore, even within the
homogeneous setup the result is stronger  -- it deals with Hausdorff
measure rather than just dimension.

\begin{theorem}\label{myth1}
Let $\vf\in C^{(n)}(I)$, $\psi$ be an approximation function and
$\lambda\in C^{(2)}(I)$. Assume that $w(f'_1,\dots,f'_n)(x)\neq 0$
for all $x\in I$. Then for any $0<s\le 1$
\begin{equation}\label{th1eq}
\HHH^s(A_n(\psi,\lambda))=\HHH^s(I)\qquad\mbox{if}\qquad
\sum_{q=1}^\infty \left(\frac{\psi(q)}{q}\right)^s\cdot q^n=\infty.
\end{equation}
\end{theorem}

Note that whenever the sum in (\ref{th1eq}) diverges, the theorem
implies that $\HHH^s(A_n(\psi,\lambda))>0$. In turn, it follows from
the definition of Hausdorff dimension that
$\dim(A_n(\psi,\lambda))\ge s$. In particular, it is easily verified
that the sum in (\ref{th1eq}) diverges whenever
$s<(n+1)/(1+\tau_{\psi})$, where
$$
\tau_{\psi}:=\liminf_{q\to\infty}\frac{-\log\psi(q)}{\log q}
$$
is the lower order of $1/\psi$ at infinity. Thus,
Theorem~\ref{myth1} readily gives the following inhomogeneous
version of the Dodson-Dickinson lower bound \cite{dod_dick} for
non-degenerate curves.

\begin{corl}\label{cor1}
Let $\vf$, $\psi$ and $\lambda$ be as in Theorem~\ref{myth1} with
$\tau_\psi\ge n$. Then
\begin{equation}\label{s}
\dim A_n(\psi,\lambda)\ge \frac{n+1}{\tau_{\psi}+1} \, .
\end{equation}
\end{corl}

Note that in the case of $\psi(q)=q^{-v}$ we have that
$\tau_{\psi}=v$ as one would expect. In the case of $\lambda$ being
a constant function and $\psi(q)=q^{-v}$, the above corollary  has
previously been obtained by the author in \cite{bod_vesti}. Bugeaud
\cite{bugeaud} has established  (\ref{s}) within the context of
approximation by algebraic integers; i.e.  in the case that $\CCC$
is the Veronese curve $(x^n,\dots,x)$ and $\lambda: x  \to x^{n+1}$.

In the case $s=1$, the $s$-dimensional Hausdorff measure $\HHH^s$ is
simply one dimensional Lebesgue measure on the real line $\RR$ .
Thus, Theorem~\ref{myth1} trivially gives rise to a complete
inhomogeneous analogue of the theorem of Beresnevich, Bernik,
Kleinbock and Margulis \cite{BBKM} in the case of non-degenerate
curves.

\begin{corl}\label{myth3}
Let $\vf$, $\psi$ and $\lambda$ be as in Theorem~\ref{myth1}.
Furthermore, suppose that the associated curve given by
(\ref{e:002}) is non-degenerate. Then
$$
|A_n(\psi,\lambda)|=|I| \quad\mbox{if}\quad\sum_{h=1}^\infty
h^{n-1}\psi(h)=\infty.
$$
\end{corl}
Here and elsewhere $|X|$ will stand for the Lebesgue measure of a
measurable subset $X$ of $\RR$.

\bigskip

\noindent\textbf{Upper bounds.} It is believed that the lower bound
for $\dim A_n(\psi,\lambda)$ given in Corollary  \ref{cor1} is
sharp. Establishing that this is the case, represents a challenging
problem and in general is open even in the homogeneous setup  -- it
has only been verified in some special cases \cite{baker,
Beresnevich-Bernik-Dodson-02:MR2069553, Bernik-1983a,
DodsonRynneVickers-1989b}. In particular,  Baker \cite{baker} has
settled the problem for planar curves within the homogeneous setup.
To the best of our knowledge,  nothing seems to be known in the
inhomogeneous case. The following result, which is an inhomogeneous
generalisation of Baker's theorem, gives a complete theory for
planar curves in the inhomogeneous case.

\begin{theorem}\label{myth2}
Let $\psi:\RR^+\to\RR^+$ be an approximation function with
$\tau_\psi\ge 2$. Let $f_1,f_2,\lambda\in C^{(2)}(I)$ be such that
the associated curve $\CCC$ given by~(\ref{e:002})${}_{n=2}$ is
non-degenerate everywhere except possibly on a  set of Hausdorff
dimension less than $\frac{3}{\tau_\psi+1}$. Then
$$
\dim A_2(\psi,\lambda)=\frac{3}{\tau_\psi+1}\,.
$$
\end{theorem}

\section{Lower bounds: Proof of Theorem~\ref{myth1}}

The proof of Theorem~\ref{myth1} will rely on the ubiquitous systems
technique as developed in \cite{ber_vel}. Essentially, the notion of
ubiquitous system  represents a convenient way of describing the
`uniform' distribution of the naturally arising points (and more
generally sets)  from a given Diophantine approximation
inequality/problem - see \cite{ber_vel, mubquity}.

\subsection{Ubiquitous systems in $\RR$}

For the sake of simplicity, we introduce a restricted notion  of
ubiquitous system, which is more than adequate for the applications
we have in mind.

\noindent Let $I_0$ be an interval in $\RR$ and $\PPP=(P_\alpha)_{\alpha\in
J}$ be a family of resonant points $P_\alpha$ of $I_0$ indexed by an
infinite set $J$. Next, let $\beta\;:\;J\to
\RR^+\;:\;\alpha\mapsto\beta_\alpha$ be a positive function on $J$.
Thus the function $\beta$ attaches a `weight' $\beta_\alpha$ to the
resonant point $P_\alpha$. Assume that for every $t\in \NN$ the set
$J_t=\{\alpha\in J: \beta_\alpha\le 2^t\}$ is finite.

Throughout, $\rho\;:\;\RR^+\to\RR^+$ will denote a function
such that $\rho(t)\to 0$ as $t\to\infty$ and it will be referred to as a {\em ubiquitous function}. Also, $B(x,r)$ will denote the
ball (or rather the interval) centered at $x$ with radius $r$.

\begin{definition}
Suppose that there exists a ubiquitous function $\rho$ and an
absolute constant $k>0$ such that for any interval $I\subseteq I_0$
$$
\liminf_{t\to\infty} \left|\bigcup_{\alpha\in J_t}B\big(P_\alpha,
\rho(2^t)\big)\cap I\right|\ge k|I|.
$$
Then the system $(\PPP, \beta)$ is called locally ubiquitous in
$I_0$ with respect to $\rho$.
\end{definition}

Let $(\PPP,\beta)$ be a ubiquitous system with respect to $\rho$ and
$\Psi$ be an approximation function. Let $\Lambda(\PPP,\beta,\Psi)$
be the set of points $\xi\in\RR$ such that the inequality
$$
|\xi-P_{\alpha}|<\Psi(\beta_\alpha)
$$
holds for infinitely many $\alpha\in J$. We will be making use of
the following lemma, which is an easy consequence of Corollary~2 (in
the case $s=1$) and Corollary~4 (in the case $s<1$) from
\cite{ber_vel}.

\begin{lemma}\label{lem1_intro}
Let $\psi$ be an approximation function and $(\PPP, \beta)$ be a
locally ubiquitous system with respect to $\rho$. Suppose there exists a real number $\lambda \in (0,1)$ such that
$\rho(2^{t+1})<\lambda \, \rho(2^t)$ for all $n\in\NN$. Then,
$$
\HHH^s(\Lambda(\PPP,\beta,\Psi))=\HHH^s(I_0) \quad\mbox{if}\quad  \sum_{t=1}^\infty\frac{\Psi(2^t)^s}{\rho(2^t)}=\infty.
$$
\end{lemma}

\subsection{Reduction of Theorem~\ref{myth1} to a ubiquity statement}

First some notation. Let $\vf=(f_1,\dots,f_n)$  be as in Theorem \ref{myth1} and denote by $\FFF_n$ the set of all functions
$$
a_0+a_1f_1(x)+a_2f_2(x)+\ldots+a_nf_n(x)
$$
where  $a_0,\ldots,a_n$ are integer coefficients, not all zero. Given
a function $F\in\FFF_n$, the {\em height} $H(F)$ of $F$ is defined as
$$H(F):=\max\{|a_1|,\ldots,|a_n|\} \ . $$  For $H>1$, let $\FFF_n(H)$ denote the
subclass of $\FFF_n$ given by
$$
\FFF_n(H)=\{F\in\FFF_n: H(F)\le H\}.
$$
Given an inhomogeneous function $\lambda$, let $R_{\lambda}=\{\alpha\in I: \exists\; F\in \mathcal{F}_n,\;
F(\alpha)+\lambda(\alpha)=0\}$. Then, for $\alpha\in R_{\lambda}$
the quantity  $H(\alpha):=\min\{H(F)\; |\; F\in \mathcal{F}_n,\;
F(\alpha)+\lambda(\alpha)=0\}$ will be referred to as the height of
$\alpha$.

\noindent To illustrate the above notions, consider the following
concrete example. Let the functions $f_i(x)=x^i$ be powers of $x$. Then
$\FFF_n$ is simply the set of all non-zero integral polynomials of
degree at most $n$. Furthermore, if $\lambda$ is identically zero,
then $R_{\lambda}$ is simply the set of algebraic numbers in $I$ of
degree at most $n$. On the other hand, if $\lambda(x)=x^{n+1}$ then $R_{\lambda}$ is
simply the set of algebraic integers in $I$ of degree exactly $n+1$.

The key to establishing Theorem~\ref{myth1} is the following ubiquity statement.

\begin{prop}\label{prop1}
The system $(R_{\lambda}, H(\alpha))$ is locally ubiquitous in $I$
with respect to $\rho(q)=q^{-n-1}$.
\end{prop}

\noindent We postpone the proof of  Proposition~\ref{prop1} to the next section. We now establish
Theorem~\ref{myth1} modulo the proposition. Note that without loss
of generality we can assume that $I$ is a closed interval. Then,
since the functions $f_i^{(j)}$ and $\lambda^{(k)},\;0\le j\le n,\;
1\le i\le n,\; 0\le k\le 2$ are continuous we have that
\begin{equation}\label{e:003}
\forall\ x\in I\qquad |f_i^{(j)}(x)|\le C,\quad
|\lambda^{(k)}(x)|\le
\end{equation}
for some absolute constant $C$. Therefore we get
$$|F'(x)|\le n\, CH(F):=MH(F).$$

\noindent Let $\alpha\in R_{\lambda}$. Then, by definition there exists a function
$F\in\mathcal{F}_n$ such that $F(\alpha)+\lambda(\alpha)=0$ and
consider the interval
$$J :=  \Big(\alpha-(2M)^{-1}H(F)^{-1}\psi(H(F)),\alpha+(2M)^{-1}H(F)^{-1}\psi(H(F))\Big).$$
For any $x\in J\cap I$, we have that
\begin{equation}\label{e:004}
 |F'(x)+\lambda'(x)|\le MH(F)+C\le 2MH(F).
\end{equation}
The latter inequality holds for all sufficiently large $H(F)$. Using
the Mean Value Theorem,  we obtain
$$
F(x)+\lambda(x)=F(\alpha)+\lambda(\alpha)+(F'(x_2)+\lambda'(x_2))(x-\alpha).
$$
By (\ref{e:004}), we get that $|F(x)+\lambda(x)|\le \psi(H(F))$ and so it follows that
whence
\begin{equation}\label{e:004a}
\Lambda(R_\lambda,H(\alpha), (2M)^{-1}H(F)^{-1}\psi(H(F)))\subset
A_n(\psi,\lambda).
\end{equation}

In view of the divergent sum condition in Theorem \ref{myth1}, we have that
$$
\sum_{t=1}^\infty
\frac{((2M)^{-1}2^{-t}\psi(2^t))^s}{2^{-t(n+1)}}\asymp
\sum_{t=1}^\infty 2^{t(n+1)}\left(\frac{\psi(2^t)}{2^t}\right)^s\ge
\sum_{h=1}^\infty h^{n}\left(\frac{\psi(h)}{h}\right)^s=\infty  .
$$
Thus,  Lemma~\ref{lem1_intro} implies that
$$ \HHH^s (\Lambda(R_\lambda,H(\alpha),
(2M)^{-1}H(F)^{-1}\psi(H(F)))) \, = \, \HHH^s(I)  \, . $$
This together with (\ref{e:004a}) implies that
$$\HHH^s(A_n(\psi,\lambda))=\HHH^s(I) \ .   $$  Modulo establishing Proposition~\ref{prop1}, this  completes the proof of
Theorem~\ref{myth1}.

\subsection{Proof of Proposition~\ref{prop1}}

Without loss of generality we can assume that $f_1(x)=x$ as
otherwise we can use the Inverse Function Theorem to change
variables and ensure the condition.
Let  $\Phi(Q,\delta)$ denote the set of $x\in I$
such that
\begin{equation}\label{e:007}
|F(x)|<\delta Q^{-n}
\end{equation}
for some $F\in \mathcal{F}_n(Q)$. We shall make use of the following lemma regarding the measure of $\Phi(Q,\delta)$.

\begin{lemma}\label{lem5}
There is an absolute constant $\delta>0$ with the following
property: for any $x_0\in I$ there is a neighborhood $I_0\subset I$
of $x_0$ such that for any interval $J\subset I_0$ there exists a
sufficiently large $Q_1>0$ such that for all $Q>Q_1$ we have
$|\Phi(Q,\delta)\cap J|<|J|/2$.
\end{lemma}

\noindent This lemma is a consequence of Theorem~2.1 in  \cite{BBKM}. Since
$I$ is compact, it is easy to see that $I_0$ can be taken to be $I$.

Take $Q_1$ and $\delta$ from Lemma~\ref{lem5}. Define
$C_1:=\delta^{\frac{1}{n+1}}$ and fix some number
$Q>\frac{1}{C_1}Q_1$. Let $\xi\in I\backslash\Phi(C_1Q,\delta)$. The
goal is to show that we can find $\alpha\in R_\lambda$
such that
\begin{equation}\label{e:008}
H(\alpha)\le K_1Q\quad\mbox{and}\quad |\xi-\alpha|\le K_2Q^{-n-1}
\end{equation}
where the constants $K_1$ and $K_2$ are independent from both $Q$
and $J$. It would immediately follow that for $Q>\frac{1}{C_1}Q_1$,
$$
\frac{|J|}{2}\le |J\backslash \Phi(C_1 Q,\delta)|\le \left|
\bigcup_{H(\alpha)\le K_1Q}B(\alpha, K_2Q^{-n-1})\cap J\right|
$$
and thus
\begin{equation}\label{vb1}
\left| \bigcup_{H(\alpha)\le Q}B(\alpha, Q^{-n-1})\cap J \right|\ge
\frac{|J|}{2K_2K_1^{n+1}}.
\end{equation}
Taking $Q=2^t$ and setting $\rho(H):=H^{-n-1}$,  inequality (\ref{vb1})
implies that $(R_{\lambda}, H(\alpha))$ is locally ubiquitous in $I$
with respect to $\rho$ -- the statement of Proposition \ref{prop1}.

We now proceed to establishing (\ref{e:008}). Consider the system of
inequalities
\begin{equation}\label{e:009}
\left\{\begin{array}{l}
|a_nf_n(\xi)+\ldots+a_1f_1(\xi)+a_0|<Q^{-n};\\[2ex]
|a_1|,|a_2|,\ldots,|a_n|\le Q.\\
\end{array}\right.
\end{equation}
It defines a convex body in $\RR^{n+1}$ symmetric about the origin.
Consider its successive minima $\tau_1,\ldots,\tau_{n+1}$. By
definition, $\tau_1\le\tau_2\le\ldots\le\tau_{n+1}$. Note that
$\tau_1>C_1$. Indeed, otherwise we would have $H\le C_1Q$ and
$$
|a_nf_n(\xi)+\ldots+a_0|\le C_1Q^{-n}= C_1^{n+1}\left(C_1
Q\right)^{-n}\le\delta H^{-n},
$$
a contradiction.
By Minkowski's theorem on successive minima \cite{lit8},
$\tau_1\cdots\tau_{n+1}\le 1$. Thus, we obtain the bound
$$
\tau_{n+1}\le (\tau_1\cdot \tau_2\cdots\tau_n)^{-1}<C_1^{-n}=C_2
$$
where $C_2$ is an absolute constant depending only on $Q$. Finally,
by the definition of $\tau_{n+1}$, we obtain  the set of $n+1$ linearly
independent functions
$F_j(X)=a_n^{(j)}f_n(X)+\ldots+a_1^{(j)}X+a_0^{(j)},\; 1\le j\le
n+1$ with integer coefficients $a_i^{(j)}$ such that
\begin{equation}\label{e:010}
\left\{\begin{array}{rcl}
|F_j(\xi)|&\le& C_2Q^{-n};\\[1ex]
|a_i^{(j)}|&\le& C_2Q,\quad i=\overline{1,n}.
\end{array}\right.
\end{equation}

\noindent Now consider the following system of linear equations
\begin{equation}\label{e:011}
\left\{\begin{array}{rcl}
\theta_1F_1(\xi)+\ldots+\theta_{n+1}F_{n+1}(\xi)+\lambda(\xi)&=&0;\\[2ex]
\theta_1F'_1(\xi)+\ldots+\theta_{n+1}F'_{n+1}(\xi)+\lambda'(\xi)&=&Q+\sum_{i=1}^n|F'_i(\xi)|;\\[2ex]
\theta_1a_j^{(1)}+\ldots+\theta_{n+1}a_j^{(n+1)}&=&0,\qquad 2\le
j\le n.
\end{array}\right.
\end{equation}

\noindent We transform this system in the following manner. Take each $j$-th row
($2\le j\le n$), multiply it by $f_j'(\xi)$ and subtract the result
from the second row. As a result, the second row will have the form
$\theta_1a_1^{(1)}+\ldots+\theta_{n+1}a_1^{(n+1)}$. Similarly we can
transform the first row to the form
$\theta_1a_0^{(1)}+\ldots+\theta_{n+1}a_0^{(n+1)}$. Since the matrix
$(a_i^{(j)}),\; 0\le i\le n,\; 1\le j\le n+1$ is non-degenerate, the
system~(\ref{e:011}) has a unique solution
$\theta_1,\ldots,\theta_{n+1}$. Choose integers $t_1,t_2,\ldots,t_n$
such that $|t_i-\theta_i|<1,\; 1\le i\le n+1$. Consider the function
$$
\begin{array}{ccc}
  F(X) &= &t_1F_1(X)+\ldots+t_{n+1}F_{n+1}(X)+\lambda(X) \\[2ex]
   & = & x_nf_n(X)+\ldots+x_1X+x_0+\lambda(X),
\end{array}
$$
where $x_i=t_1a_i^{(1)}+\ldots+t_{n+1}a_i^{(n+1)}$.
By the first equation in (\ref{e:011}), we obtain
$$|F(\xi)|<(n+1)C_2Q^{-n}=C_3Q^{-n} \ . $$ Further, by the second equation
in (\ref{e:011}), we obtain $|F'(\xi)|>Q$ and
$$
\begin{array}{ll}
|F'(\xi)|<Q+2\sum_{i=1}^{n+1}|F'_i(\xi)|&\le
Q+2(n+1)((n+1)C_2\cdot CQ)\\
&=(1+2(n+1)^2 C_2\cdot C)Q=C_4Q.
\end{array}
$$
Now consider the coefficients $x_i$. They are obviously integers. By
the third equation in (\ref{e:011}), we get $|x_m|\le (n+1)C_2
Q=C_5Q$ for all $m\ge 2$. The bounds for $x_0$ and $x_1$ are given
by
$$
\begin{array}{rcl}
|x_1| &\le& |F'(\xi)|+|\lambda(\xi)|+\sum_{i=2}^n|x_if'_i(\xi)|\\[2ex]
 & \le & C_4Q+(n-1)(n+1)C_2Q\cdot C+C\le C_6Q
\end{array}
$$
and
$$
\begin{array}{rcl}
|x_0| & \le & |F(\xi)|+|\lambda(\xi)|+\sum_{i=1}^n|x_if_i(\xi)|\\[2ex]
 & \le & C_3Q^{-n}+(n-1)(n+1)|C_2\cdot CQ+C_6CQ+C\le C_7Q
\end{array}
$$
for sufficiently large $Q$. Thus, for every $\xi\in I\backslash\Phi(C_1Q,\delta)$ there exists
$F(x)\in\mathcal{F}_n$ such that
\begin{equation}\label{e:012}
\left\{\begin{array}{l} |F(\xi)+\lambda(\xi)|\le C_3Q^{-n};\\[2ex]
Q\le |F'(\xi)+\lambda'(\xi)|\le C_4Q;\\[2ex]
|H(F)|\le \max(C_5,C_6,C_7)Q.
\end{array}\right.
\end{equation}
It is easy to check that $\max\{C_5,C_6,C_7\}=C_7$. Hence,
$|H(F)|\le C_7Q$ or equivalently $F\in \mathcal{F}_n(C_7Q)$.

The next goal  is to show that the function
$F(x)+\lambda(x)$ constructed above has a root $\alpha$ satisfying
conditions \eqref{e:008}.

\begin{lemma}\label{lem6}
Let $\sigma(F)$ be a set of all $x\in I$ satisfying the following
system of inequalities:
$$
\left\{\begin{array}{l} |F(x)+\lambda(x)|\le C_3Q^{-n}\\[1ex]
Q\le|F'(x)+\lambda'(x)|\le C_4Q,
\end{array}\right.
$$
where $F\in \mathcal{F}_n(C_7Q)$. Let $Q$ satisfy the condition
$$
(n\cdot C\cdot C_7Q+C)\cdot 2C_3Q^{-n-1}+C\le \frac{1}{2}Q.
$$ 
Then for all $x_0\in\sigma(F)\cap [a+2C_3Q^{-n-1},b-2C_3Q^{-n-1}]$
there exists a number $\alpha\in(x_0-2C_3Q^{-n-1},x_0+2C_3Q^{-n-1})$
such that $F(\alpha)+\lambda(\alpha)=0$.
\end{lemma}

\noindent\textit{Proof.}  By the Mean Value Theorem,
$$
F'(x)+\lambda'(x)=F'(x_0)+\lambda'(x_0)+(F''(x_1)+\lambda''(x_1))(x-x_0),
$$
where $x_1$ is some point between $x$ and $x_0$. Taking
$x\in(x_0-2C_3Q^{-n-1},x_0+2C_3Q^{-n-1})$ and using (\ref{e:003}) we
get $F''(x_1)+\lambda''(x_1)\le n\cdot C\cdot C_7Q+C$. Therefore
$$|(F''(x_1)+\lambda''(x_1))(x-x_0)|\le (n\cdot C\cdot C_7Q+C)\cdot
2C_3Q^{-n-1}\le \frac{1}{2}Q-C.$$ Finally we get that for all real
$x$ such that $|x-x_0|\le 2C_3Q^{-n-1}$ the following inequality is
satisfied
$$
\begin{array}{rcl}
|F'(x)+\lambda'(x)|&\ge&
|F'(x_0)|-|\lambda'(x_0)|-|(F''(x_1)+\lambda''(x))(x-x_0)|\\[2ex]
&>&|F'(x_0)|/2.
\end{array}
$$
In particular it means that the function $F'(x)+\lambda'(x)$ has the
same sign within the given interval. Again, on  using the Mean Value
Theorem we get that
$F(x)+\lambda(x)=F(x_0)+\lambda(x_0)+(F'(x_2)+\lambda'(x_2))(x-x_0)$,
where $x_2$ lies between $x$ and $x_0$. Set $x=x_0\pm 2C_3q^{-n-1}$.
Then
$$
\begin{array}{rcl}
|(F'(x_2)+\lambda'(x_2))(x-x_0)| & > & 2C_3Q^{-n-1}|F'(x_0)|/2\\[2ex]
 & \ge & C_3Q^{-n}\ge |F(x_0)+\lambda(x_0)|.
\end{array}
$$
Note that for the two different values of $x$ the expression
$$
(F'(x_2)+\lambda'(x_2))\cdot(x-x_0)
$$
has different signs. Therefore
the value of
$F(x)+\lambda(x)=F(x_0)+\lambda(x_0)+(F'(x_2)+\lambda'(x_2))(x-x_0)$
has different signs at the two ends of the interval
$$
[x_0-2C_3Q^{-n-1},x_0+2C_3Q^{-n-1}].
$$
Thus the function $F(x)+\lambda(x)$ has a root within this interval and thereby  completes the proof of Lemma~\ref{lem6}.

\hfill\hfill$\boxtimes$

\

In view of Lemma \ref{lem6},  we have that for all $\xi$ satisfying
system (\ref{e:012}) there exists $\alpha$ with $H(\alpha)\le C_7Q$
such that
$$
F(\alpha)+\lambda(\alpha)=0
$$
and
$$ |\xi-\alpha|<2C_3Q^{-n-1}.
$$

\noindent Finally, for all $\xi\in I\backslash\Phi(C_1 Q,\delta)$ we have
constructed a function $F(x)\in\FFF_n$ such that \eqref{e:012} is
satisfied. Therefore, by taking $K_1=C_7$ and $K_2=2C_3$, we find a
number $\alpha\in R_\lambda$ satisfying \eqref{e:008}. This
completes the proof of the Proposition~\ref{prop1}.


\section{Upper bounds: Proof of Theorem~\ref{myth2}}

\subsection{Preliminary notes}

First of all note that, by Corollary~\ref{cor1}, it suffices to
establish the lower bound
\begin{equation}\label{vb8c}
 \dim A_2(\psi,\lambda)\le \frac{3}{\tau_\psi+1}\,.
\end{equation}
Note that there is nothing to prove if $\tau_\psi=2$. Thus, without
loss of generality we can assume that $\tau_\psi>2$. Further, the
definition of $\tau_\psi$ readily implies that for any $v<\tau_\psi$
we have that $\psi(q)\ll q^{-v}$ for all sufficiently large $q$. It
follows that for any $v<\tau_\psi$ we have that
$A_2(\psi,\lambda)\subset A_2(v,\lambda)$. Therefore, (\ref{vb8c})
will follow if we consider the special case of $\psi(q)=q^{-v}$ with
$2<v<\tau_\psi$ and let $v\to\tau_\psi$. Therefore, from now on we
fix a $v>2$ and concentrate on establishing the bound
\begin{equation}\label{vb8d}
 \dim A_2(v,\lambda)\le \frac{3}{v+1}\,.
\end{equation}

\subsection{Auxiliary lemmas}

As in the proof of Proposition \ref{prop1}, there is no loss of generality in  assuming
that $f_1(x)=x$. Then we simply denote $f_2(x)$ by $f(x)$. With the
aim of establishing Theorem~\ref{myth2} we fix $v>2$. By the
conditions of Theorem~\ref{myth2}, we have that $f''(x)\not=0$ for
all $x$ except a set of Hausdorff dimension $\le\frac{3}{v+1}$.
Using the standard arguments -- see \cite{vbb} -- we can assume
without loss of generality that
\begin{equation}\label{vb8a}
    c_1\le |f''(x)|\le c_2\qquad\text{for all }x\in I,
\end{equation}
where $c_1,c_2$ are positive constants.

\begin{lemma}[Pyartly \cite{pyartly}]\label{lem3}
Let $\delta,\nu>0$ and $I\subset \RR$ be some interval. Let
$\phi(x)\in C^n(I)$ be a function such that $|\phi^{(n)}(x)|>\delta$
for all $x\in I$. Then there exists a constant $c(n)$ which depends
only on $n$, such that
$$
|\{x\in I\; :\; |\phi(x)|<\nu\}|\le
c(n)\left(\frac{\nu}{\delta}\right)^{\frac{1}{n}}.
$$
\end{lemma}

Before stating the next lemma recall that $\FFF_2$ is the set of all
functions of the form $a_0+a_1x+a_2f(x)$, where $a_0,a_1,a_2$ are
integers not all zero; $H=H(F)=\max\{|a_1|$, $|a_2|\}$.

\begin{lemma}\label{lem4}
There are constants $C_1>0$ and $\epsilon_0>0$ such that  for all
$F\in\FFF_2$ and any subinterval $J\subset I$ of length
$|J|\le\epsilon_0$ at least one of the following inequalities is
satisfied for all $x\in J$:
$$|F'(x)+\lambda'(x)|>C_1H(F)\quad\mbox{or}\quad |F''(x)+\lambda''(x)|>C_1H(F).$$
\end{lemma}

\noindent\textit{Proof.}  For the case of $\lambda(x)\equiv 0$ this
is proved in \cite[Lemmas 5, 6]{vbb}. To finish the proof in
inhomogeneous case it is sufficient to note that $|\lambda'(x)|\ll
1$ and $|\lambda''(x)|\ll 1$.\\  \hspace*{\fill}$\boxtimes$

In what follows without loss of generality we can assume that
$|I|\le\epsilon_0$ -- see \cite{vbb} for analogues arguments.

\begin{lemma}\label{lem9}
Fix some $0\le\delta\le 1$ and a positive number $H$. Denote by
$N(\delta)$ the number of triples $(a_0,a_1,a_2)\in \ZZ^3$
satisfying $\max\{|a_i|:i=0,1,2\}\le H$ such that there exists a
solution $x\in I$ to the system
\begin{equation}\label{e:013}
\left\{
\begin{array}{l}
|F(x)+\lambda(x)|\le H^{-v}\\[1ex]
|F'(x)+\lambda'(x)|\le H^\delta.
\end{array}\right.
\end{equation}
Then for $v>0$, $N(\delta)\ll H^{1+\delta}$.
\end{lemma}

\noindent\textit{Proof.}  Since $|\lambda'(x)|\ll 1$,
$|\lambda(x)|\ll 1$ and $\delta\ge 0$, we have that (\ref{e:013})
implies the following system
\begin{equation}\label{e:014}
\left\{
\begin{array}{l}
|F(x)|\ll H^\delta\\[1ex]
|F'(x)|\ll H^\delta.
\end{array}\right.
\end{equation}

\noindent Subtracting the second inequality of (\ref{e:014}) multiplied by $x$
from the first inequality of (\ref{e:014}) gives
\begin{equation}\label{vb8b}
\left\{
\begin{array}{l}
|a_0+a_2(f(x)-xf'(x))|\ll H^\delta\\[1ex]
|a_1+a_2f'(x)|\ll H^\delta.
\end{array}\right.
\end{equation}

If $|a_1|=H$ then we have $2H+1$ possibilities for $a_2$. By
(\ref{vb8a}), for each fixed pair $(a_1,a_2)$ the interval of $x$
satisfying the second inequality of (\ref{vb8b}) is of length
$O(H^\delta a_2^{-1})$. Therefore, $a_2(f(x)-xf'(x))$ may vary on an
interval of length $O(H^\delta)$ only. Hence, for every fixed pair
$(a_1,a_2)$ we have $O(H^\delta)$ possibilities for $a_0$. Thus, we
have $O(H^{\delta+1})$ triples $(a_0,a_1,a_2)$ with $|a_1|=H$.

Consider the case $|a_2|=H$. Note that by (\ref{vb8a}), $f'(x)$ is
strictly monotonic and finite. Therefore one can change variables by
setting $t=-f'(x)$; $f(x)-xf'(x)=g(t)$. Note that, by (\ref{vb8a}),
the variable $t$ belongs to some finite interval $J$. Furthermore,
the function $g(t)$ is bounded, continuously differentiable on $J$
and $|g'(t)|\ll 1$. Therefore, the system (\ref{vb8b}) transforms to
$$
\left\{
\begin{array}{l}
\left|a_0+a_2g(t)\right|\ll H^\delta;\\[3pt]
\left|a_1-a_2t\right|\ll H^\delta;\\[3pt]
t\in J.
\end{array}\right. \Longrightarrow
\left\{
\begin{array}{l}
\left|\frac{a_0}{a_2}+g(t)\right|\ll H^{\delta-1};\\[7pt]
\left|\frac{a_1}{a_2}-t\right|\ll H^{\delta-1};\\[5pt]
t\in J.
\end{array}\right.
$$

\noindent Note that
$$g\left(\frac{a_1}{a_2}\right) =g(t+\Delta)=g(t)+\Delta
g'(\xi)=g(t)+O(H^{\delta-1})  \ , $$ where $\Delta=\frac{a_1}{a_2}-t$.
Hence all solutions of the system are also solutions of the
inequality
$$
\left|g\left(\frac{a_1}{a_2}\right)+\frac{a_0}{a_2}\right|\ll
H^{\delta-1}.
$$

\noindent One can easily check that for $|a_2|=H$ the number of integer
solutions of this inequality is not greater than $CH^{1+\delta}$ for
some constant $C$. Therefore $N(\delta)\ll H^{1+\delta}$ and the
proof is complete. \hfill$\boxtimes$

\begin{lemma}\label{lem12}
Consider the plane defined by the equation $Ax+By+Cz=D$ where
$A,B,C,D$ are integers with $(A,B,C)=1$. Then the area $S$ of any
triangle on this plane with integer vertices is at least $\frac12
\sqrt{A^2+B^2+C^2}$.
\end{lemma}

\noindent\textit{Proof.} Denote by ${\bf x}, {\bf y}$ and ${\bf z}$
some points on the considered plane not all lying on the same line.
Take one more integer point ${\bf v}$ somewhere outside the plane.
We now calculate the volume $V$ of the tetrahedron ${\bf xyzv}$.

On one hand the volume of every tetrahedron with integer vertexes is
at least $\frac16$. Therefore $V\ge \frac16$.

On the other hand, $V=\frac13Sh$, where $S$ is the area of the
triangle ${\bf xyz}$ and $h$ is the distance between ${\bf v}$ and
the plane. Therefore,
$$
\frac16\le V=\frac13 Sh \Longleftrightarrow \frac2h\le S.
$$
Let ${\bf v}=(\alpha,\beta,\gamma)$. Then
$$
h=\frac{|A\alpha+B\beta+C\gamma-D|}{\sqrt{A^2+B^2+C^2}}\ge
\frac{1}{\sqrt{A^2+B^2+C^2}},
$$
since $\alpha, \beta$ and $\gamma$ are integers and $h>0$. Thus,
$S\ge \frac12\sqrt{A^2+B^2+C^2}$ as required. \hfill$\boxtimes$


\subsection{Proof of Theorem \ref{myth2}}

Let $\sigma:=\frac{3}{v+1}$ be the required bound in (\ref{vb8d}).
The strategy of the proof is to construct a collection of coverings
$D_i=\{d_{ij}:{j\in J}\}$ of $A_2(v,\lambda)$ by intervals $d_{ij}$
such that for any $\epsilon>0$
$$
\sum_{j\in J}|d_{ij}|^{\sigma+\epsilon}\to 0\quad\mbox{as}\quad i\to
\infty.
$$
The bound (\ref{vb8d}) will then follow from the definition of
Hausdorff dimension. Note that $A_2(v,\lambda)$ can be represented
in one of the following forms
\begin{equation}\label{e:015}
A_2(v,\lambda)=\bigcap_{n=1}^\infty\bigcup_{H=n}^\infty
A(a_0,a_1,a_2)\quad\mbox{ and}
\end{equation}
$$
A_2(v,\lambda)=\bigcap_{n=1}^\infty\bigcup_{t=n}^\infty B(t),
$$
where $A(a_0,a_1,a_2)$ is the set of $x\in I$ satisfying
\begin{equation}\label{eqt}
|a_0+a_1x+a_2f(x)+\lambda(x)|<H^{-v}
\end{equation}
for the particular triple $(a_0,a_1,a_2)$;
$$B(t)=\bigcup\limits_{2^{t-1}\le H<2^t} A(a_0,a_1,a_2);\quad
H=\max\{|a_1|,|a_2|\}.
$$
Therefore for any $n\in \NN$ the collection of sets $A(a_0,a_1,a_2)$
with $H\ge n$ is a covering of $A_2(v,\lambda)$. Analogously for any
$n\in \NN$ the collection of $B(t)$ with $t\ge n$ is a covering of
$A_2(v,\lambda)$.

Fix some positive small number $\epsilon$. Divide every set
$A(a_0,a_1,a_2)$ into three subsets:
\begin{equation}\label{e:016}
A_1(a_0,a_1,a_2)=\{x\in A(a_0,a_1,a_2):
|F'(x)+\lambda'(x)|>H^{1-\epsilon}\};
\end{equation}
\begin{equation}\label{e:017}
A_2(a_0,a_1,a_2)=\left\{x\in A(a_0,a_1,a_2):
H^{\frac{2-v}{3}}<|F'(x)+\lambda'(x)|\le H^{1-\epsilon} \right\};
\end{equation}
\begin{equation}\label{e:018}
A_3(a_0,a_1,a_2)=\left\{x\in A(a_0,a_1,a_2): |F'(x)+\lambda'(x)|\le
H^{\frac{2-v}{3}} \right\}.
\end{equation}

For any of these collections we can construct the associated sets
$A_2^{(1)}(v,\lambda)$, $A_2^{(2)}(v,\lambda)$ and
$A_2^{(3)}(v,\lambda)$ analogously to $A_2(v,\lambda)$ -- see
\eqref{e:015}. One can easily check that
$$
A_2(v,\lambda)=A_2^{(1)}(v,\lambda)\cup A_2^{(2)}(v,\lambda)\cup
A_2^{(3)}(v,\lambda).
$$
Therefore it is sufficient to prove (\ref{vb8d}) for
$A_2(v,\lambda)$ replaced by either of these subsets.

\vspace{2ex}\textbf{The set $A_2^{(1)}(v,\lambda)$.} Since
$|\lambda(x)|\ll1$, we have that
$$
|a_1+a_2f'(x)+\lambda(x)|>H^{1-\epsilon} \qquad\Longrightarrow\qquad
|a_1+a_2f'(x)|\gg H^{1-\epsilon}.
$$
Since $|f''(x)|>d$ for all $x\in I$, we have that $a_1+a_2f'(x)$ is
a monotonic function. Therefore the set of $x\in I$ such that
$|a_1+a_2f'(x)|>H^{1-\epsilon}$ is a union of at most two intervals.
For one interval we have that
$$a_1+a_2f'(x)\ll -H^{1-\epsilon}$$
and for the other we have that
$$a_1+a_2f'(x)\gg H^{1-\epsilon}.$$
We see that the sign of $F'(x)+\lambda'(x)$ on each of these
intervals doesn't change. Therefore $F(x)+\lambda(x)$ is monotonic
on them, where $F(x)=a_0+a_1x+a_2f(x)$.

Thus for sufficiently large $H$ the set $A_1(a_0,a_1,a_2)$ is a
union of at most 2 intervals (note that it can be empty, i.e. be a
union of empty intervals).

Using Lemma~\ref{lem3} and inequality in (\ref{e:016}) we get that
the length of each interval is $\ll H^{-v-1+\epsilon}$.

We will use the following cover of $A_2^{(1)}(v,\lambda)$:
$$
C_n=\bigcup_{H=n}^\infty A_1(a_0,a_1,a_2).
$$
Note that for a fixed $H$ the number of different pairs $(a_1,a_2)$
is no greater than $4H$. By \eqref{eqt} there are $O(H)$
possibilities for $a_0$ if $(a_1,a_2)$ are fixed. Therefore an
appropriate $s$-volume sum for $C_n$ will be

$$
C\ll \sum_{H=n}^\infty H^2\cdot
H^{s(\epsilon-1-v)}=\sum_{H=n}^\infty H^{2-s(1+v-\epsilon)}.
$$
This sum tends to zero as $n\to\infty$ in case of
$2-s(1+v-\epsilon)<-1$, that is $s>\frac{3}{1+v-\epsilon}$. Thus,
\begin{equation}\label{e:019}
\dim(A_2^{(1)}(v,\lambda))\le \frac{3}{1+v-\epsilon}.
\end{equation}

\vspace{2ex}\textbf{The set $A_2^{(2)}(v,\lambda)$.} Here we have
the inequality $|F'(x)+\lambda'(x)|\le H^{1-\epsilon}$. Therefore
Lemma~\ref{lem4} implies
\begin{equation}\label{e:020}
\forall x\in A_2(a_0,a_1,a_2),\quad |F''(x)+\lambda''(x)|\gg H.
\end{equation}
In other words $|a_2f''(x)+\lambda''(x)|\gg H$. This implies that
$|a_2|\gg H$. Note that~\eqref{e:020} is also true in the case of
$x\in A_3(a_0,a_1,a_2)$.

Let $\delta$ be an arbitrary number in $(0,1]$. Consider the set
\begin{equation}\label{e:021}
A_\delta(v,\lambda)=\bigcap_{n=1}^\infty\bigcup_{H=n}^\infty
A_\delta(a_0,a_1,a_2),
\end{equation}
where $A_{\delta}(a_0,a_1,a_2)$ is the set of $x\in
A_2(a_0,a_1,a_2)$ with the following property
\begin{equation}\label{e:022}
H^{1-\frac{1}{3}(v+1)\delta}<|F'(x)+\lambda'(x)|\le H^{1-\delta}.
\end{equation}

We have that $|F''(x)+\lambda''(x)|\gg H$. Therefore, the set
$A_\delta(a_0,a_1,a_2)$ consists of at most 4 intervals. Consider
the following cover $C_n$ for $A_\delta(v,\lambda)$:
$$
C_n=\bigcup_{H=n}^\infty A_\delta(a_0,a_1,a_2).
$$

By Lemma~\ref{lem9}, for a fixed $H$ there exist only
$O(H^{2-\delta})$ nonempty sets $A_{\delta}(a_0,$ $a_1, a_2)$. By
Lemma \ref{lem3} the length of each interval in
$A_\delta(a_0,a_1,a_2)$ is bounded by
$H^{-v-1+\frac{1}{3}(v+1)\delta}$. Therefore the corresponding
$s$-volume sum for $C_n$ is bounded by
\begin{equation}\label{vb8f}
\sum_{H=n}^\infty H^{2-\delta}\cdot
H^{s\left(-v-1+\frac{1}{3}(v+1)\delta\right)}.
\end{equation}
If $s>\frac{3}{v+1}$ then the exponent of $H$ is equal to
$$
2-\delta+s\left(-v-1+\frac{1}{3}(v+1)\delta\right)<2-\delta-3+\delta=-1.
$$
Hence for $s>\frac{3}{v+1}$ the right hand side of (\ref{vb8f})
tends to 0 as $n\to\infty$. It follows that
$\dim(A_\delta(v,\lambda))\le \frac{3}{v+1}$ for any $\delta \in
[0,1]$.

For simplicity denote by $k$ the quantity $\frac{1}{3}(v+1)$. Note
that $k>1$. For $l\in\NN$ consider the set $A_2^{(2)}(v,\lambda)$ as
a union
\begin{equation}\label{vb8g}
A_2^{(2)}(v,\lambda)=\bigcup_{i=1}^l A_{\delta_i}(v,\lambda)\cup
A_{\delta^*}(v,\lambda),
\end{equation}
where $\delta_1=\epsilon, \delta_{i+1}=k\delta_i, \delta^*=1$. Since
$k>1$ we have that $\delta_i\to \infty$ as $i\to\infty$. Therefore
there exists a natural number $l$ which depends on $\epsilon$ only
such that $\delta_{l+1}>1$ and $\delta_l\le 1$. This proves
(\ref{vb8g}).

Since the Hausdorff dimension of each set $A_{\delta_i}(v,\lambda)$
and $A_{\delta^*}(v,\lambda)$ appearing in (\ref{vb8g}) is not
greater than $\frac{3}{v+1}$, we get that
$\dim(A_2^{(2)}(v,\lambda))\le\frac{3}{v+1}$.

\vspace{2ex}\textbf{The set $A_2^{(3)}(v,\lambda)$.} Consider the
set
$$
B_3(t):=\bigcup_{2^{t-1}\le H<2^t}A_3(a_0,a_1,a_2)
$$
and let $\delta:=\frac{2-v}{3}$.

Recall that for all $x\in A_3(a_0,a_1,a_2)$ we have
$|F''(x)+\lambda''(x)|\gg H$. Therefore, by Lemma \ref{lem3}, we
have that the length of each interval in $A_3(a_0,a_1,a_2)$ with
$H\asymp 2^t$ is not greater than
$$
(H^{-v}/H)^{\frac{1}{2}}\ll 2^{-t\left(\frac{v+1}{2}\right)}.
$$

Fix a sufficiently small positive number $\epsilon_1$. Let
$c=1+\epsilon_1$. For every $t$ divide the interval $I$ into
$2^{ct}$ equal subintervals of length $2^{-ct}|I|\ll 2^{-ct}$. These
subintervals are divided into two classes:
\begin{itemize}
\item Class I intervals. They include at most
$O\left(2^{t(\frac{3}{2}-c)}\right)$ segments from $B_3(t)$.
\item Class II intervals. They include those which are not in class I.
\end{itemize}
According to this classification consider the sets
$$
A_\mathrm{I}(v,\lambda)=\bigcap_{n=1}^\infty\ \bigcup_{t=n}^\infty \
\bigcup_{\scriptscriptstyle\text{ class I intervals $J$}} B_3(t)\cap
J;
$$
$$
A_{\mathrm{II}}(v,\lambda)=\bigcap_{n=1}^\infty\
\bigcup_{t=n}^\infty \ \bigcup_{\scriptscriptstyle\text{ class II
intervals $J$}} B_3(t)\cap J.
$$

It follow that
$$
A_2^{(3)}(v,\lambda)=A_\mathrm{I}(v,\lambda)\cup
A_\mathrm{II}(v,\lambda).
$$
The required upper bound for $A_\mathrm{I}(v,\lambda)$ will follow
on showing the following lemma.

\begin{lemma}\label{lem10}
$\dim(A_\mathrm{I}(v,\lambda))\le\frac{3}{v+1}(1+\epsilon_1)$.
\end{lemma}

\noindent\textit{Proof.}  Consider a class I interval $J$. We have
at most $O\left(2^{t(\frac{3}{2}-c)}\right)$ segments from $B_3(t)$
on it. Therefore there are not greater than $O(2^{\frac{3}{2}t})$
intervals from $B_3(t)$ lying inside class I intervals. Consider the
following cover of $A_{\mathrm{I}}(v,\lambda)$:
$$
C_n:=\bigcup_{t=n}^\infty \ \bigcup_{\scriptstyle\text{class I
intervals } J} B_3(t)\cap J.
$$
Its $\frac{3}{v+1}(1+\epsilon_1)$-volume is bounded by
$$
\sum_{t=n}^\infty 2^{\frac{3}{2}t}\cdot
2^{-t\left(\frac{v+1}{2}\right)\cdot
\frac{3(\epsilon_1+1)}{v+1}}=\sum_{t=n}^\infty
2^{\left(\frac{3}{2}-\frac{3}{2}(\epsilon_1+1)\right)t}=\sum_{t=n}^\infty
2^{-\frac{3}{2}\epsilon_1t}.
$$
It obviously tends to zero as $n\to\infty$. This finishes the proof
of the lemma.

\hfill$\boxtimes$

\bigskip

Let $J$ be a class II interval and
$F(x)=a_0+a_1x+a_2f(x)\in\mathcal{F}_2$ with $2^{t-1}\le H(F)<2^t$
and $A_3(a_0,a_1,a_2)\cap J\not=\emptyset$. Then
$$
|F(x_0)+\lambda(x_0)|\ll 2^{-vt}\qquad\text{and}\qquad
|F'(x_0)+\lambda'(x_0)|\ll 2^{\delta t}
$$
for some $F\in \FFF_2$ and $x_0\in J$. Then for all $x\in J$ we have
$$
|F'(x)+\lambda'(x)|=|(F'+\lambda')(x_0)+(x-x_0)(F''+\lambda'')(\xi)|\ll
2^{\delta t}+2^{(1-c)t}.
$$
$$
\begin{array}{l@{}l}|F(x)\!+\!\lambda(x)|&=|(F\!+\!\lambda)(x_0)+(x-x_0)(F'+\lambda')(x_0)+(x-x_0)^2(F''+\lambda'')(\xi)|\\[2ex]
&\ll 2^{-v t}+2^{(\delta-c)t}+2^{(1-2c)t}.
\end{array}
$$
Choose a sufficiently small $\epsilon_1>0$ such that
\begin{equation}\label{e:023}
v>2+3\epsilon_1.
\end{equation}
Then we have
$$
2^{\delta
t}<2^{(1-c)t}\qquad\text{and}\qquad2^{(\delta-c)t}<2^{(1-2c)t}.
$$
One can see that $2^{-vt}$ is always less than the other summands
$2^{(\delta-c)t}$ and $2^{(1-2c)t}$. Hence in the case of
(\ref{e:023}) we get the inequalities
\begin{equation}\label{e:024}
|F(x)+\lambda(x)|\ll 2^{(1-2c)t},
\end{equation}
\begin{equation}\label{e:025}
|F'(x)+\lambda'(x)|\ll 2^{(1-c)t}
\end{equation}
for all $x\in J$.

\begin{lemma}\label{lem11}
For every fixed $J$ as above all points
$\vec{a}=(a_0,a_1,a_2)\in\ZZ^3$ such that $A_3(a_0,a_1,a_2)\cap
J\neq \emptyset$ lie on a single affine plane.
\end{lemma}

\noindent\textit{Proof.}  Suppose that there exist four integer
points $\vec{a},\vec{b},\vec{c},\vec{d}$ not lying on the same plane
such that $A_3(\vec{a})\cap J\neq\emptyset,A_3(\vec{b})\cap
J\neq\emptyset, A_3(\vec{c})\cap J\neq\emptyset$ and
$A_3(\vec{d})\cap J\neq\emptyset$. It means that the points
$\vec{a},\vec{b},\vec{c},\vec{d}$ form a tetrahedron with integer
vertexes. Therefore its volume is at least $\frac{1}{6}$.

On the other hand all of these four points must lie inside a
parallelepiped~$R$ formed by the inequalities (\ref{e:024}),
(\ref{e:025}) and $|a_2|<H$ for a fixed $x\in J$. The volume of this
figure is bounded by
$$
V\ll 2\cdot2^{t(1-2c)}\cdot2\cdot2^{t(1-c)}\cdot2\cdot2^t\cdot
D^{-1}\ll 2^{t(3-3c)}\cdot D^{-1},
$$
where $D$ is the determinant of the matrix
$$
\left(\begin{array}{ccc} 1&x&f(x)\\
0&1&f'(x)\\
0&0&1
\end{array}\right)
$$
i.e. $D=1$. Since $c>1$ we have $V=o(1)$ contrary to $V\ge1/6$. The
proof is complete.

\hfill$\boxtimes$

Let the plane from Lemma \ref{lem11} have the form $Ax+By+Cz=D$. We
evaluate the intersection area of this plane with parallelepiped $R$
specified in the proof of the lemma. In order to do this let us
consider the body $P_\Delta$ given by the inequalities
\begin{equation}\label{e:026}
\left\{\begin{array}{l}
|F(x)+\lambda(x)|\le 2^{t(1-2c)};\\
|a_2|\le 2^t;\\
|Aa_0+Ba_1+Ca_2-D|\le \Delta,
\end{array}\right.
\end{equation}
where $\Delta>0$ is a positive parameter. Here $a_0,a_1,a_2$ are
viewed as real variables. The volume of $P_{\Delta}$ can be
expressed in two different ways. Firstly, since the determinant of
system (\ref{e:026}) is $B-Ax$, we have that
\begin{equation}\label{e:027}
V(P_\Delta)=\frac{8\cdot 2^{t(2-2c)}\cdot \Delta}{|B-Ax|} .
\end{equation}
Secondly, $V(P_\Delta)=S\cdot h$ where $S$ is the area of the edges
defined by the third inequality of (\ref{e:026}) and $h$ is a
distance between these edges. That is
\begin{equation}\label{e:028}
V(P_\Delta)=S\cdot \frac{2\Delta}{\sqrt{A^2+B^2+C^2}}.
\end{equation}
Hence on combining \eqref{e:027} and \eqref{e:028} we obtain that
\begin{equation}\label{e:029}
S\asymp \frac{2^{t(2-2c)}\cdot\sqrt{A^2+B^2+C^2}}{|B-Ax|}.
\end{equation}
Note that $S$ is the area of the intersection of the plane
$Aa_0+Ba_1+Ca_2-D=0$ with the figure defined by the first two
inequalities of (\ref{e:026}). Therefore the intersection area of
this plane with the parallelepiped is not greater than $S$. Note
that all points $\va$ should lie inside this intersection and
\eqref{e:029} gives an estimate for its area.

\

\noindent\textbf{Case (i):} We consider intervals $J$ of type II
such that not all points $\va$ associated with $J$ lie on the same
line. By Lemma \ref{lem12}, we get that the number of such points on
a fixed interval $J$ is bounded by
\begin{equation}\label{e:030}
N\ll
\frac{2^{t(2-2c)}\cdot\sqrt{A^2+B^2+C^2}}{|B-Ax|}/\sqrt{A^2+B^2+C^2}=
\frac{2^{t(2-2c)}}{|B-Ax|}.
\end{equation}
Since $J$ is not a class I interval we get that for all $x\in J$,
\begin{equation}\label{e:031}
|B-Ax|\ll 2^{t(\frac{1}{2}-c)}.
\end{equation}

Similarly to (\ref{e:026}) we consider two more systems of
inequalities:
$$
\left\{\begin{array}{l}
|F(x)+\lambda(x)|\le 2^{t(1-2c)};\\
|Aa_0+Ba_1+Ca_2-D|\le \Delta;\\
|F'(x)+\lambda'(x)|\le 2^{t(1-c)}\\
\end{array}\right.\mbox{ and }
\left\{\begin{array}{l}
|F'(x)+\lambda'(x)|\le 2^{t(1-c)};\\
|a_2|\le 2^t;\\
|Aa_0+Ba_1+Ca_2-D|\le \Delta.\\
\end{array}\right.
$$
Analogously we get additional bounds for $N$, namely
\begin{equation}\label{e:032}
 N\ll
\frac{2^{t(2-3c)}}{|T|}\qquad\text{and}\qquad
N\ll\frac{2^{t(2-c)}}{|A|},
\end{equation}
where
$$
T=\det\left(\begin{array}{ccc} 1&x&f(x)\\
A&B&C\\
0&1&f'(x)\\
\end{array}\right)=f'(x)(B-Ax)-(C-Af(x)).
$$
Since $J$ is a class II interval then
$$
|f'(x)(B-Ax)-(C-Af(x))|\ll 2^{t(\frac{1}{2}-2c)}.
$$
This result with (\ref{e:031}) implies
\begin{equation}\label{e:033}
|C-Af(x)|\ll 2^{t(\frac{1}{2}-c)}.
\end{equation}

The second inequality in (\ref{e:032}) together with the fact that
$J$ is a class II interval implies $|A|\ll 2^{\frac{1}{2}t}$.

Fix $A$. Denote by $M(A)$ the number of possible integer triples
$(A,B,C)$ which can be the coefficients of a plane corresponding to
some class II interval. It follows from (\ref{e:031}) and
(\ref{e:033}) that
$$
M(A)\le |I|\cdot |A|+1\ll |A|.
$$
In fact it is the number of fractions $\frac{B}{A}$ in the interval
$I$. Parameter $C$ is uniquely defined by $A$ and $B$.

Suppose there exist two class II intervals $J_1$ and $J_2$ with the
same coefficients $(A,B,C)$ of appropriate plane. Applying
inequality (\ref{e:031}) we get
$$
\left\{\begin{array}{l} \forall x\in J_1, |B-Ax|\ll
2^{t(\frac{1}{2}-c)};\\
\forall y\in J_2, |B-Ay|\ll 2^{t(\frac{1}{2}-c)}.
\end{array}\right. \Rightarrow |A(x-y)|\ll 2^{t(\frac{1}{2}-c)}
\Rightarrow |x-y|\ll \frac{2^{t(\frac{1}{2}-c)}}{|A|}.
$$
Therefore for a fixed $(A,B,C)$ the number $x$ can lie only in the
interval of the length $\frac{2^{t(\frac{1}{2}-c)}}{|A|}$. Finally
we get that the number of class II intervals associated with the
triple $(A,B,C)$ is at most $\frac{2^{\frac{1}{2}t}}{|A|}$.

We will use the following cover for the $A_\mathrm{II}(v,\lambda)$:
$$
C_n=\bigcup_{t=n}^\infty \bigcup_{J \mbox{\tiny\:are class II
intervals}} B_3(t)\cap J;
$$
Using the second inequality in (\ref{e:032}) to estimate the number
of intervals in $B_3(t)\cap J$ we get that the $s$-volume sum for
this cover is bounded by
$$
C=\sum_{t=n}^\infty \sum_{(A,B,C)} \frac{2^{\frac{1}{2}t}}{|A|}\cdot
\frac{2^{t(2-c)}}{|A|}\cdot 2^{-st(\frac{v+1}{2})}=\sum_{t=n}^\infty
2^{t(\frac{3}{2}-\epsilon_1-s(\frac{v+1}{2}))}\sum_{(A,B,C)}
\frac{1}{|A|^2},
$$
where $(A,B,C)$ run through possible coefficients of planes
corresponding to type II intervals under consideration. Transforming
this series we get
\begin{equation}\label{e:034}
\begin{array}{l@{}l}\displaystyle\sum_{t=n}^\infty 2^{t(\frac{3}{2}-\epsilon_1-s(\frac{v+1}{2}))}
\sum_{|A|=1}^{2^{\frac{t}{2}}}\frac{1}{|A|}&\displaystyle\ll
\sum_{t=n}^\infty t\cdot
2^{t(\frac{3}{2}-\epsilon_1-s(\frac{v+1}{2}))}\\
&\displaystyle\ll \sum_{t=n}^\infty
2^{t(\frac{3}{2}-s(\frac{v+1}{2}))}.
\end{array}
\end{equation}

If $s>\frac{3}{v+1}$ then this series obviously tends to zero as
$n\to \infty$.

\

\noindent\textbf{Case (ii):} We consider intervals $J$ of type II
such that all points $\va$ associated with $J$ lie on the same line
$L$.  Fix such an interval $J$. Represent this line in a form:
$$
\mathbf{\alpha}+ t \mathbf{\beta}
$$
where $\mathbf{\alpha}=(\alpha_0,\alpha_1,\alpha_2)$ is an integer
point on $L$, $\mathbf{\beta}=(\beta_0,\beta_1\,beta_2)$ is a vector
between the nearest integer points on $L$ and $t$ is an arbitrary
real number. Then all the vectors $(a_0,a_1,a_2)$ associated with
$J$ are of the form
$$
a_0=\alpha_0+k\beta_0, a_1=\alpha_1+k\beta_1, a_2=\alpha_2+k\beta_2,
$$
where $\alpha_i,\beta_i$ are fixed and $k\in \ZZ$ vary. Since $J$ is
of class II, there are at least $2^{t(\frac{3}{2}-c)}$ different
values $k$. For each vector $(a_0,a_1,a_2)$ under consideration we
have that $|a_0|\ll 2^t$. Hence taking values of $|a_0|$ for two
different vectors for $J$ and subtracting one value from another we
get
\begin{equation}\label{e:035}
|\beta_0(k_1-k_2)|\ll 2^t \Rightarrow |\beta_0|\le
2^{t(c-\frac{1}{2})}.
\end{equation}
Similarly we obtain the same inequalities for $\beta_1$ and
$\beta_2$.

Now consider inequalities (\ref{e:024}) and (\ref{e:025}) for the
same two vectors. Again subtracting one inequality from the other we
get
\begin{equation}\label{e:036}
\left\{
\begin{array}{l}
|(k_1-k_2)(\beta_0+\beta_1x+\beta_2f(x))|\le 2\cdot
2^{t(1-2c)};\\[2ex]
|(k_1-k_2)(\beta_1+\beta_2f'(x))|\le 2\cdot 2^{t(1-c)}.
\end{array}\right.
\end{equation}

Since there are at least $2^{t(\frac32-c)}$ different values $k$, we
can ensure that $|k_1-k_2|\ge 2^{t(\frac32-c)}$ for some $k_1, k_2$.
Dividing these inequalities by $|k_1-k_2|$ and changing variables as
in Lemma~\ref{lem9} give the system
$$ \left\{
\begin{array}{l}
|\beta_0+\beta_2g(y)|\le 2\cdot 2^{-\frac{t}{2}};\\[1ex]
|\beta_1+\beta_2y|\le 2\cdot 2^{-\frac{t}{2}}
\end{array}\right.
$$
where $y=f'(x)$ and $g(y)=f(x)-xf'(x)$. Using (\ref{e:035}) we get
$H=\max\{|\beta_0|$, $|\beta_1|, |\beta_2|\}\ll
2^{t(c-\frac{1}{2})}$. Substituting this into the system we get
\begin{equation}\label{e:037}
\left\{
\begin{array}{l}
|\beta_0+\beta_2g(y)|\ll
H^{\frac{1}{1-2c}}=H^{-\frac{1}{1+2\epsilon_1}};\\[2ex]
|\beta_1+\beta_2y|\ll
H^{\frac{1}{1-2c}}=H^{-\frac{1}{1+2\epsilon_1}}.
\end{array}\right.
\end{equation}
For a fixed value of $\beta_2$ the number of possibilities for
$\beta_1$ is $O(|\beta_2|)$. For a fixed $\beta_1$ and $\beta_2$ the
number of possibilities for $\beta_0$ is $O(1)$. Denote by
$K(\beta_2)$ the number of solutions $(\beta_0,\beta_1,\beta_2)$ of
(\ref{e:037}), where $\beta_2$ is fixed. Then we get that
$K(\beta_2)\ll |\beta_2|$.

%
%
%

Suppose that for two different intervals $J_1$ and $J_2$ the
parameters $\beta_0,\beta_1$ and $\beta_2$ coincide. Using the
second inequality in (\ref{e:037}) we get
$$
|\beta_2(y_1-y_2)|\ll H^{-\frac{1}{1+2\epsilon_1}}
$$
where $y_1\in f'(J_1)$ and $y_2\in f'(J_2)$. Since for all $x\in I,
f'(x)>d>0$ the inequality can be transformed to the form
$$
|x_1-x_2|\ll \frac{H^{-\frac{1}{1+2\epsilon_1}}}{|\beta_2|}\ll
\frac{2^{-\frac{1}{2}t}}{|\beta_2|}.
$$
Therefore the number of class $B$ intervals $J$ with parameters
$\beta_0,\beta_1,\beta_2$ is not greater than
\begin{equation}\label{e:038}
\frac{2^{t(c-\frac{1}{2})}}{|\beta_2|}.
\end{equation}
Further, since $|\alpha_2+k\beta_2|\le 2^t$ there are at most
$\frac{2^t}{|\beta_2|}$ intervals inside $J\cap B_3(t)$. Using
(\ref{e:038}) we have the following upper bound for $s$-volume sum.
$$
\begin{array}{ll}C&\displaystyle=\sum_{t=n}^\infty \sum_{(\beta_0,\beta_1,\beta_2)}
\frac{2^{t(c-\frac{1}{2})}}{|\beta_2|}\cdot
\frac{2^t}{|\beta_2|}\cdot 2^{-st(\frac{v+1}{2})}=\sum_{t=n}^\infty
2^{t(\frac{3}{2}+\epsilon_1-s(\frac{v+1}{2}))}\sum_{(\beta_0,\beta_1,\beta_2)}
\frac{1}{|\beta_2|^2}\\[15pt]
&\displaystyle \ll\sum_{t=n}^\infty
2^{t(\frac{3}{2}+\epsilon_1-s(\frac{v+1}{2}))}\sum_{|\beta_2|\le
2^{t(c-\frac12)}}\frac{K(\beta_2)}{|\beta_2|^2}\\
&\displaystyle \ll \sum_{t=n}^\infty t\cdot
2^{t(\frac{3}{2}+\epsilon_1-s(\frac{v+1}{2}))}.
\end{array}
$$
Using the same arguments as in the case when $(a_0,a_1,a_2)$ lie on
a plane we obtain that
$$
C\ll \sum_{t=n}^\infty
2^{t(\frac{3}{2}+2\epsilon_1-s(\frac{v+1}{2}))}.
$$

Combining this series with (\ref{e:034}) we get an estimate
$$
\dim(A_\mathrm{II}(v,\lambda))\le \frac{3+4\epsilon_1}{v+1}.
$$

Therefore finally for $v>2+3\epsilon_1$ we get that
$$
\begin{array}{l@{}l}
\dim (A_2(v,\lambda))&=\dim(A_2^{(1)}(v,\lambda)\cup
A_2^{(2)}(v,\lambda)\cup A_\mathrm{I}(v,\lambda)\cup
A_\mathrm{II}(v,\lambda))\\[2mm]
&\displaystyle \le\max\left\{\frac{3}{1+v-\epsilon},
\frac{3+4\epsilon_1}{v+1},\frac{3}{v+1}(1+\epsilon_1)\right\}.
\end{array}
$$

Since $\epsilon$ and $\epsilon_1$ can be made arbitrary small then
all values in the maximum can be made arbitrary close to
$\frac{3}{v+1}$, thus implying (\ref{vb8c}) and thereby completing
the proof of Theorem~\ref{myth2}.

\

\noindent\textit{Acknowledgements. } DB is grateful to Victor
Beresnevich, Vasili Bernik and Sanju Velani for introducing me to
the wonderland of metrical Diophantine approximation  and for their
numerous helpful discussions.

{\footnotesize

\

\sc \noindent Dzmitry Badziahin\\
\rm Mathematics Department, University of York, Heslington, York,
YO10
5DD, England\\
\it Email address: \tt db528@york.ac.uk

}

\end{document}